\documentclass[conference]{IEEEtran}
\IEEEoverridecommandlockouts
\usepackage{cite}
\usepackage{amsmath,amssymb,amsfonts}
\usepackage{algorithmic}
\usepackage{algorithm}
\usepackage{graphicx}
\usepackage{textcomp}
\usepackage{multicol}
\usepackage{xcolor}
\usepackage{url}
\def\BibTeX{{\rm B\kern-.05em{\sc i\kern-.025em b}\kern-.08em
    T\kern-.1667em\lower.7ex\hbox{E}\kern-.125emX}}
\begin{document}

\title{Bilevel Optimization for Smart Irrigation Control: A Framework for Enhancing Water-Use Efficiency in Agriculture\\
}

\author{\IEEEauthorblockN{Roman Chertovskih}
\IEEEauthorblockA{\textit{SYSTEC-ARISE} \\
\textit{Faculty of Engineering,}\\
\textit{University of Porto}\\
 Porto, Portugal \\
roman@fe.up.pt}
\and
\IEEEauthorblockN{ Nathalie Khalil}
\IEEEauthorblockA{\textit{SYSTEC-ARISE} \\
\textit{Faculty of Engineering,}\\
\textit{University of Porto}\\
 Porto, Portugal \\
nathalie@fe.up.pt}
\and
\IEEEauthorblockN{ Luís Sousa}
\IEEEauthorblockA{\textit{SYSTEC-ARISE} \\
\textit{Faculty of Engineering,}\\
\textit{University of Porto}\\
 Porto, Portugal \\
luismsousa103@gmail.com}

}

\maketitle

\begin{abstract}
Global warming has intensified water scarcity, posing a critical challenge for the agricultural sector, where seasonal demand often leads to significant wastage. This work addresses the need for efficient water management by proposing a smart irrigation control system based on a bilevel optimization algorithm. The framework aims to minimize water consumption across multiple fields simultaneously or, under supply constraints, to distribute available water equitably to minimize deviations from required levels.
To achieve this, the problem is modeled hierarchically: the upper level acts as a central controller managing daily water limits and field-wise allocation, while the lower level acts as a centralized scheduler determining the optimal cooperative irrigation execution across all fields. This structure ensures that both global resource constraints and local irrigation needs are met efficiently. Preliminary results from digital simulation tools suggest that the proposed framework significantly improves water-use efficiency and supports sustainable irrigation practices, particularly in water-constrained scenarios.
\end{abstract}


\section{Introduction}
Smart agriculture is an emerging frontier in sustainable resource management, aiming to optimize crop yields while minimizing environmental impact. Among the critical challenges in this domain is managing water -- an increasingly scarce resource that requires precise delivery to maintain soil health and crop productivity. While traditional irrigation systems often rely on uniform application or simple rule-based triggers, modern precision irrigation systems offer a promising alternative by accounting for spatial variability in soil types, topography, and crop needs.

Despite the potential of precision irrigation, significant hurdles remain in large-scale deployment. Key challenges include managing water distribution across heterogeneous field sectors, respecting total system capacity constraints, and ensuring autonomous operation under weather conditions. These difficulties are magnified in multi-point irrigation configurations, where coordinating water delivery across $N$ distinct sectors introduces complex trade-offs between local moisture needs and global resource limits. Furthermore, the dynamic nature of soil-water interactions, consisting of infiltration, evapotranspiration, and runoff, requires continuous, real-time optimization to maximize water-use efficiency.

To address these spatial and resource-based complexities, this work proposes the Bilevel Multi-Point Irrigation Model. This model adopts a hierarchical centralized control architecture that separates macro-resource planning from micro-operational scheduling. To formalize this interaction, the system utilizes a bi-level optimization framework where a central coordinator guides a global operational scheduler. While bilevel programming is increasingly used in areas such as supply chain \cite{avraamidou2017}, network design~\cite{madadi2023}, telecommunication \cite{nguyen2026}, energy markets \cite{moghadam2022}, defense systems \cite{yuan2025}, among others, its application to irrigation systems, characterized by multi-point spatial dynamics and nonlinear soil physics, represents a novel approach to agricultural decision support.

\section{Bilevel Optimization: Overview}

\subsection{The Problem}
Bilevel optimization is formally defined as a hierarchical mathematical program \cite{sinha2017reviewbilevel, dempe2020-book, zemkoho2021}, where one optimization task is nested within another. This structure involves an upper-level decision vector $x_u \in X_U \subseteq \mathbb{R}^n$ controlled by a leader and a lower-level decision vector $x_l \in X_L \subseteq \mathbb{R}^m$ controlled by a follower. The leader seeks to optimize an upper-level objective function $F: X_U \times X_L \rightarrow \mathbb{R}$ while accounting for the follower's rational response. Mathematically, the problem is formulated as follows:

\begin{equation}
\begin{aligned}
\min_{x_u \in X_U, x_l \in X_L} \quad & F(x_u, x_l) \\
\text{subject to} \quad & G_k(x_u, x_l) \le 0, \quad k = 1, \dots, K \\
& x_l \in \Psi(x_u).
\end{aligned}
\end{equation}

For every decision $x_u$ made by the leader, the follower reacts by solving a parametric optimization problem to find an optimal response. This response is represented by the set-valued mapping $\Psi: X_U \rightrightarrows X_L$, which defines the lower-level reaction set as:

\begin{align}\nonumber
\Psi(x_u) = & \arg\min_{x_l \in X_L} \{ f(x_u, x_l) : g_j(x_u, x_l) \le 0, \\ & \quad j = 1, \dots, J \}\label{reaction set}
\end{align}
where $f: X_U \times X_L \rightarrow \mathbb{R}$ and $g_j$ denote the follower's objective and constraints, respectively. A solution pair $(x_u, x_l^*)$ is considered feasible for the leader only if $x_l^* \in \Psi(x_u)$.

\subsection{Solution Types}
A defining feature of bilevel problems is their inherent asymmetry. While both levels have distinct objectives and constraints, the leader typically has complete knowledge of the follower's problem, enabling anticipation of responses. Conversely, the follower only observes the leader’s decision before optimizing their own strategy. Depending on the assumptions regarding the follower's cooperation, two primary formulations exist: 

\begin{itemize}
    \item \textbf{Optimistic Position:} The leader assumes the follower will choose the response from $\Psi(x_u)$ that best serves the upper-level objective.
    \item \textbf{Pessimistic Position:} The leader optimizes for the worst-case scenario, assuming the follower may choose a response from $\Psi(x_u)$ that is most detrimental to the upper-level goal.
\end{itemize}
In this paper, we will adopt an optimistic approach.

\section{Our Irrigation Framework}
The core innovation of the proposed model lies in applying the hierarchy described above to an irrigation context. In this ``Leader-Follower'' structure, the Upper Level (Leader) chooses strategic variables that shape the operational context, while the Lower Level (Follower) optimizes local performance in response. 

At the Strategic Level (Upper Level), the system acts as a coordinator, determining the optimal dynamic soil moisture safety thresholds ($X_{\min, j, i}$) for each sector over time. This level balances the global objective of maintaining field capacity against the hard constraints of total water availability ($C_{\text{total}}$).At the Operational Level (Lower Level), the system acts as a centralized scheduler that solves a joint, multi-field Optimal Control Problem (OCP). This operational layer minimizes the total collective water consumption across all sectors simultaneously, ensuring that local soil conservation dynamics and safety margins dictated by the Upper Level are met globally. By treating the safety threshold as a strategic, time-varying variable rather than a fixed parameter, the system shifts from a rigid controller to an adaptive decision-support architecture.

\section{Problem Formulation}

\subsection{Upper Level: Strategic Resource Coordinator}
The goal of the central coordinator is to manage macro-level water allocations during supply deficits. It minimizes the aggregate squared deviation of the dynamically adjusted reference thresholds $X_{\min, j, i}$ from the ideal agronomic crop water requirements ($\text{need}_j$) across all sectors and time steps.

\begin{equation}\label{obj:deviation}
\min_{X_{\min}, \text{deviation}} \quad J = \sum_{j=1}^{N} \sum_{i=1}^{M} (\text{deviation}_{i,j})^2
\end{equation}
subject to:
\begin{align}
& \sum_{j=1}^{N} u_{i,j} \le C_{\text{total}}, \quad \forall i=1, \dots, M \label{eq:total_cap} \\
& 0 \le X_{\min, j, i} \le \text{need}_{j}, \quad \forall j, i \label{eq:domain_bound}\\
& \text{deviation}_{i,j} \ge \text{need}_{j} - X_{\min, j, i}, \quad \forall j, i \label{eq:slack_bound}\\
& \text{deviation}_{i,j} \ge 0, \quad \forall j, i \\
& (u,x) \in \Psi(X_{\min}). \label{eq:lower_link}
\end{align}

The coordinator adjusts the dynamic reference matrix $X_{\min, j, i} \in \mathbb{R}^{N \times M}$, which defines the minimum allowable moisture baseline for field $j$ at step $i$. Constraint (\ref{eq:total_cap}) restricts aggregate daily irrigation to the delivery infrastructure capacity $C_{\text{total}}$. Equations (\ref{eq:domain_bound})--(\ref{eq:slack_bound}) establish a penalization mechanism: when shared water is abundant, $X_{\min, j, i}$ equals $\text{need}_j$ and deviation drops to zero. During severe droughts, the leader strategically reduces $X_{\min, j, i}$ to enable controlled deficit irrigation, thereby minimizing the systemic stress penalty quadratically.

\subsection{Lower Level: Operational Optimal Control}
The lower level acts as the localized irrigation execution agent. For a given threshold baseline profile $X_{\min, j, i}$ dictated by the leader, it minimizes total volumetric water applications while preserving local soil dynamics and respecting structural capacity boundaries.

\begin{equation}
\min_{u,x, \text{loss}, \text{overflow}} \quad \sum_{j=1}^{N} \sum_{i=1}^{M} u_{i,j}
\end{equation}
subject to:
\begin{align}
& x_{i,j} \ge X_{\min, j, i}, \quad \forall i, j \label{eq:state_constraint} \\
& x_{1,j} = x_{S,j}, \quad \forall j \\
& x_{i+1,j} = x_{i,j} + h f(t_i, x_{i,j}, u_{i,j}, x_{i+1,j}), \quad \forall i=1, \dots, M-1 \label{eq:dynamics} \\
& \text{overflow}_{i,j} \ge (x_{i,j} - x_{FC})c_{s,j}, \quad \forall i, j \label{eq:overflow_constraint} \\
& \text{loss}_{i,j} \ge k \cdot x_{i,j}, \quad \forall j, i = 1, \dots, M-1 \label{eq:loss_constraint} \\
& u_{i,j} \ge 0, \quad x_{i,j} \ge 0, \quad \text{overflow}_{i,j} \ge 0, \quad \text{loss}_{i,j} \ge 0
\end{align}
The set-valued mapping $\Psi(X_{\min})$ represents the collection of lower-level rational responses passed back up to couple the hierarchical system layers.

Equation (\ref{eq:dynamics}) represents the soil moisture dynamics, which are governed by the function $f$ that integrates irrigation, crop evapotranspiration, and environmental factors, defined as:

\begin{equation}
\begin{aligned}
& f(t_i, x_{i,j}, u_{i,j}, x_{i+1,j}) =  K_{I} u_{i,j} [1 - R_{F} \sin(\alpha_j)]\\& - K_{C} \text{evtp}_{0}(t_{i}) 
 + [1 - R_{F} \sin(\alpha_j)] \text{Rainfall}(t_{i}) \\
& -\text{loss}_{i,j} - \text{overflow}
_{i,j}
\end{aligned}
\end{equation}

The dynamics are composed of four primary terms:
\begin{itemize}
    \item \textbf{Irrigation and Rainfall:} Both water inputs are adjusted by the factor $[1 - R_{F} \sin(\alpha_j)]$ to account for runoff losses determined by the soil slope $\alpha_j$.
    \item \textbf{Crop Requirements:} The term $K_{C} \text{evtp}_{0}(t_{i})$ represents the volume of water extracted by the crop through evapotranspiration.
    \item \textbf{Soil Losses and Overflow:} The functions $\text{loss}_{i,j}$ and $\text{overflow}_{i,j}$ account for deep infiltration when $x \le x_{FC}$ and saturation-related losses when $x > x_{FC}$.
    \item $K_I, K_C$ are coefficients for irrigation type and crop needs, while $R_F$ is a runoff coefficient which quantifies the fraction of water that fails to infiltrate the soil and instead flows over the land surface.
\end{itemize}

\begin{table*}[tp]
\centering
\caption{Indices, Variables, and Parameters}
\begin{tabular}{lll}
\hline
\textbf{Symbol} & \textbf{Description} & \textbf{Unit} \\ \hline
\textit{Indices} & & \\
$i$ & Discrete time steps ($1, \dots, M$) & hr or day \\
$j$ & Spatial irrigation fields/sectors ($1, \dots, N$) & -- \\
\textit{Decision Variables} & & \\
$X_{\min, j, i}$ & Dynamic strategic soil moisture threshold for sector $j$ at time $i$ & mm \\
$u_{i,j}$ & Irrigation applied at sector $j$ during time $i$ & mm \\
$x_{i,j}$ & Actual soil moisture state at sector $j$ during time $i$ & mm \\
$\text{deviation}_{i,j}$ & Upper-level water deficit tracking slack variable for sector $j$ at time $i$ & mm \\
$\text{overflow}_{i,j}$ & Lower-level saturation runoff overflow tracking variable for sector $j$ at time $i$ & mm \\
$\text{loss}_{i,j}$ & Lower-level deep percolation soil water loss variable for sector $j$ at time $i$ & mm \\
\textit{Parameters} & & \\
$\text{need}_j$ & Ideal agronomic soil moisture target/requirement for crop in sector $j$ & mm \\
$x_{FC}$ & Soil moisture threshold at maximum Field Capacity & mm \\
$x_{S,j}$ & Initial soil moisture from sensors at sector $j$  & mm \\
$C_{\text{total}}$ & Total shared water delivery capacity limit across all sectors  & mm \\
$\alpha_j$ & Slope angle of field sector $j$ & deg \\
$R_F$ & Runoff coefficient  & -- \\
$c_{s,j}$ & Excess water overflow loss efficiency coefficient for sector $j$ & -- \\
$K_I, K_C$ & Coefficients for irrigation application type and specific crop needs  & -- \\
$k$ & Constant empirical coefficient related to deep soil permeability & -- \\ \hline
\end{tabular}
\end{table*}

\section{Distinction from Existing Optimal Control Formulations}
Our proposed bilevel multi-point irrigation model extends the single-level irrigation problem presented in \cite{pereira2022irrigation} by shifting the functional role of the safety threshold from a static parameter to a dynamic strategic decision variable $X_{\min, j, i}$. In the existing literature, $x_{\min}$ is treated as a fixed value representing rigid hydric requirements for crop survival, and the problem is solved as a single-level optimal control problem. This approach required the local system to satisfy a pre-determined moisture target, defined by the constraint $x_{i} \ge x_{\min}$ (cf. \cite{pereira2022irrigation}), regardless of global resource availability. Such a configuration can lead to mathematical infeasibility or total system failure during severe droughts when the available water supply cannot meet fixed survival thresholds. 

To overcome such scenarios, our work introduces a hierarchical coordination layer that, at the upper coordinator level, optimizes the strategic dynamic soil moisture threshold matrix $X_{\min, j, i}$ for each field sector. This allows the coordinator to dynamically adjust safety margins across $N$ field sectors, facilitating controlled deficit irrigation scenarios and the spatial redistribution of water based on real-time constraints such as total system capacity $C_{\text{total}}$ or varying soil slopes $\alpha_{j}$. Consequently, while the irrigation problem in \cite{pereira2022irrigation} focuses on minimizing localized water use at a single point, our bilevel approach transforms the system into an adaptive decision-support architecture capable of balancing competing needs across an entire multi-point farm configuration.

\section{Computational platform and numerical methods}
Numerical solution of the optimal control problem was carried out in Julia~\cite{bezanson2017julia} using JuMP~\cite{dunning2017jump} as the algebraic modeling layer and IPOPT~\cite{wachter2006ipopt} as the nonlinear programming solver. JuMP provides a high-level, symbolic interface for formulating constrained optimization problems in a syntax close to mathematical notation, while delegating low-level communication with the underlying optimizer to the optimizer itself. A critical advantage of JuMP in the context of optimal control is its integration with automatic differentiation, which provides exact first- and second-order derivative information (Jacobians and Hessians) required by gradient-based solvers without manual differentiation.

Julia was chosen because it combines high productivity with performance suitable for scientific computing, and it is widely recommended for numerical work in applied mathematics. Its just-in-time (JIT) compilation yields execution speeds comparable to statically compiled languages such as C and Fortran, which is highly advantageous for the iterative nature of numerical optimization algorithms. 

The problem was formulated as a large-scale nonlinear programming problem, a standard and effective approach for optimal control when control and state variables are represented on a finite mesh. The resulting nonlinear program was solved with IPOPT, an interior-point method designed for large-scale nonlinear optimization. IPOPT uses a primal-dual barrier framework, a filter line-search strategy, second-order corrections, and restoration phases to improve robustness and convergence for constrained nonlinear problems. This methodology is particularly well-suited to optimal control problems after direct transcription, where sparse Jacobian and Hessian structures can be exploited efficiently. The algorithm implemented in IPOPT was introduced by Wächter and Biegler~\cite{wachter2006ipopt} and remains a standard reference for large-scale nonlinear programming in scientific computing.

All computations were executed on one x86 compute node of the Deucalion supercomputer \cite{deucalion}, hosted by the Minho Advanced Computing Center (MACC) in Guimarães, Portugal. The Deucalion is a heterogeneous system with separate ARM, x86, and GPU partitions; the x86 partition consists of 500 nodes, each equipped with 2 AMD EPYC 7742 processors, 256 GB DDR4 memory, and 100 Gb/s NVIDIA ConnectX-6 network connectivity. The system runs Rocky Linux 8 and uses Slurm for resource management and job scheduling. Its storage and network architecture are designed for high-throughput scientific workloads, with a Lustre parallel file system and dedicated high-speed interconnects. 

The combination of the highly optimized JuMP/IPOPT software stack and the dense multi-core, high-memory capacity of the AMD EPYC architecture within the Deucalion cluster ensured that the primary computational bottlenecks, most notably the resolution of sparse linear systems required at each interior-point iteration, were handled with maximum efficiency and numerical stability. This computational setup is well aligned with reproducible scientific workflows: the mathematical model is described at a high level in JuMP, the solver is open-source and extensively documented, and the execution platform is a modern HPC resource with clearly specified hardware and software characteristics.

\section{Algorithm}
The mathematical problem is implemented using BilevelJuMP \cite{bilevelJump2024}, an extension of the JuMP modeling language designed specifically for solving bilevel problems. The algorithm is described in detail below.  

\begin{algorithm*}[t!]
\caption{Bilevel Hierarchical Allocation for Smart Field Irrigation}
\label{alg:bilevel_irrigation}

\begin{multicols}{2}
\begin{algorithmic}[1]

\renewcommand{\algorithmicrequire}{\textbf{Input:}}
\renewcommand{\algorithmicensure}{\textbf{Output:}}

\REQUIRE ~\\
Time horizon $M$, Spatial fields $N$\\
Requirements: $\mathbf{\text{need}_j} \in \mathbb{R}^N$\\
Capacity: $C_{\text{total}} \in \mathbb{R}$\\
Initial moisture: $\mathbf{x_{S,j}} \in \mathbb{R}^N$\\
Profiles: $\mathbf{\text{rain}}, \mathbf{\text{wind}}, \mathbf{\text{temp}} \in \mathbb{R}^M$\\
Physics: $\mathbf{K_I}, \mathbf{K_C}, \mathbf{\alpha_j}, R_F, x_{FC}, k$

\ENSURE ~\\
Strategic matrix: $X_{\min} \in \mathbb{R}^{N \times M}$\\
Centralized scheduling matrix: $u \in \mathbb{R}^{M \times N}$\\
State trajectory matrix: $x \in \mathbb{R}^{M \times N}$

\STATE \COMMENT{Phase 1: Environmental Pre-processing}
\FOR{$i = 1$ \TO $M$}
    \STATE Calculate reference evapotranspiration $\text{ET}_{0,i}$ via FAO-56.
\ENDFOR

\STATE \COMMENT{Phase 2: Monolithic Bilevel Optimization Setup}
\STATE Initialize Bilevel Container via \texttt{BilevelJuMP}.

\STATE \textbf{Begin Upper-Level Problem (Central Coordinator):}
\STATE Define strategic variables: $X_{\min, j, i}, \text{deviation}_{i,j}$
\STATE Set Upper Objective: Minimize aggregate baseline deviations \eqref{obj:deviation}
\STATE Enforce Global Constraints:
\begin{itemize}
    \item Shared delivery network capacity limits \eqref{eq:total_cap}
    \item Strategic domain bounds and slacks \eqref{eq:domain_bound}--\eqref{eq:slack_bound}
\end{itemize}

\vfill\null
\columnbreak

\STATE \quad \textbf{Begin Linked Lower-Level Problem (Central Scheduler):}
\STATE \quad Define operational variables: $u_{i,j}, x_{i,j}, \text{loss}_{i,j}, \text{overflow}_{i,j}$
\STATE \quad Set Lower Objective: Minimize total collective water use
\STATE \quad Enforce Joint Multi-Field Dynamics ($\forall j \in 1\dots N, \forall i \in 1\dots M$):
\begin{itemize}
    \item Bind states to strategic matrix: $x_{i,j} \ge X_{\min, j, i}$
    \item Initialize baseline moisture states: $x_{1,j} = x_{S,j}$
    \item Evaluate mass-balance soil state transitions \eqref{eq:dynamics}
    \item Constrain deep percolation and saturation slacks
\end{itemize}
\STATE \quad \textbf{End Lower-Level Problem}

\STATE \textbf{End Upper-Level Problem}

\STATE \COMMENT{Phase 3: Telemetry \& Execution}
\STATE Set underlying solver to \texttt{Ipopt} and invoke optimization.
\IF{Solver Status is Optimal}
    \STATE Extract unified decision matrices $X_{\min}^*$, $u^*$, $x^*$
    \STATE Compute macro water-saving and yield-stress metrics.
    \STATE \RETURN $X_{\min}^*, u^*, x^*$
\ELSE
    \STATE \RETURN \texttt{Optimization Failure Flag}
\ENDIF

\end{algorithmic}
\end{multicols}
\end{algorithm*}

\section{Results}
\begin{figure}[htbp]
    \centering
    \includegraphics[width=0.45\textwidth, keepaspectratio]{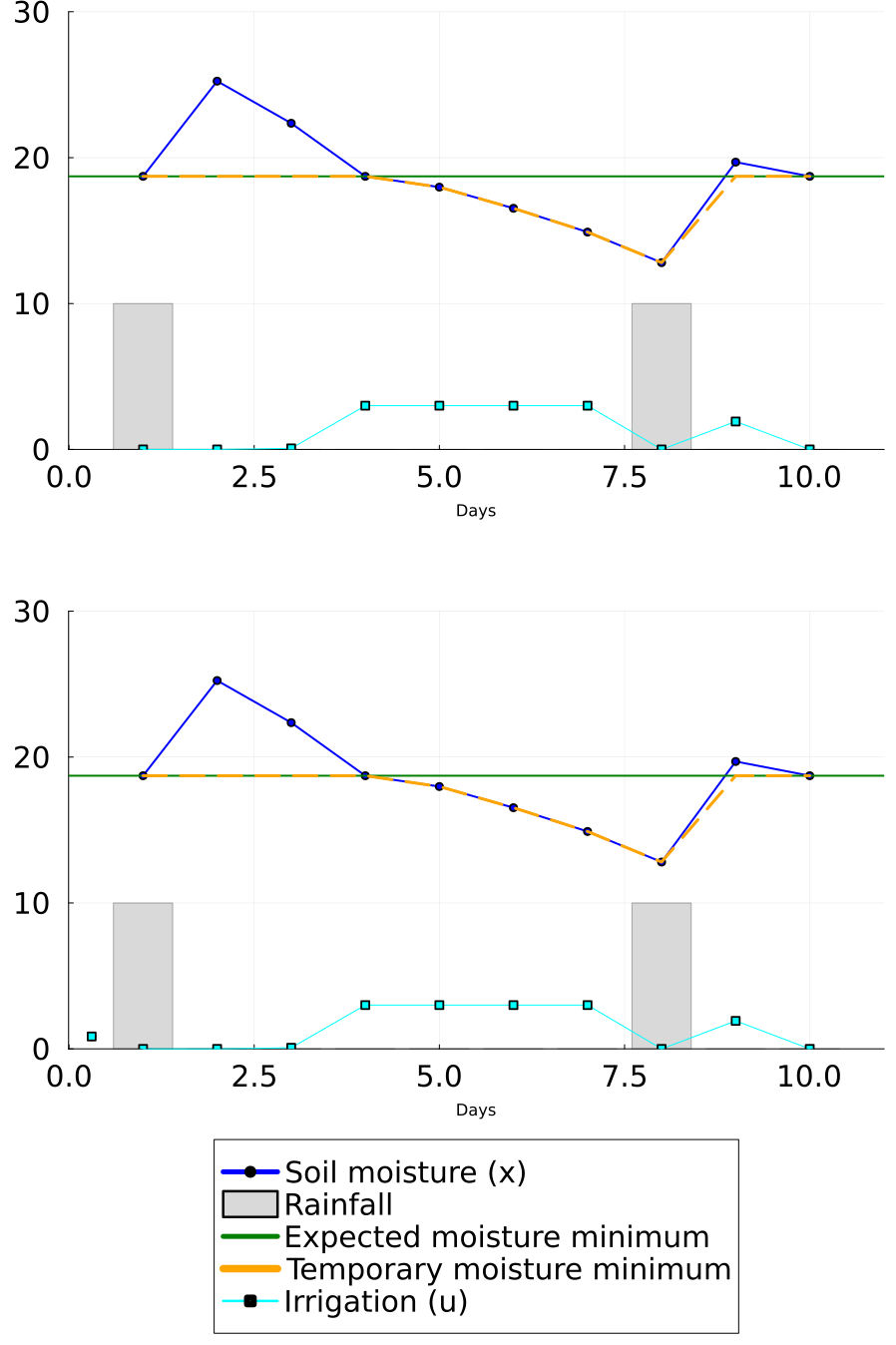}
    \caption{Irrigation profile for two identical fields under shared capacity constraints ($C_{\text{total}} = 6$~mm).}
    \label{fig:irrigation1}
\end{figure}

First, the simulation was tested across two fields with identical parameter values (soil slope, temperature, wind, etc.) under a restricted total shared daily water delivery limit $C_{\text{total}} = 6$~mm. The optimization trajectory is shown in Fig.~\ref{fig:irrigation1}.
The green line illustrates the ideal baseline requirement ($\text{need}_j$) to preserve full plant health. The dashed orange line shows the dynamically tracking temporary lower threshold ($X_{\min, j, i}$), lowered by the leader to cope with water availability restrictions. The dark blue plot tracks actual soil moisture ($x_{i,j}$), the cyan line captures irrigation applications ($u_{i,j}$), and the gray bars illustrate natural rainfall inputs. Because initial allocations and structural configurations match exactly across both sectors, their resulting schedules are identical. Following rainfall on day 1, moisture levels are naturally preserved, delaying scheduled controller actions until day 3 to minimize water shortages across the shared network. Indeed, irrigation starts at day 3.

\begin{figure}[htbp]
    \centering
    \includegraphics[width=0.45\textwidth, keepaspectratio]{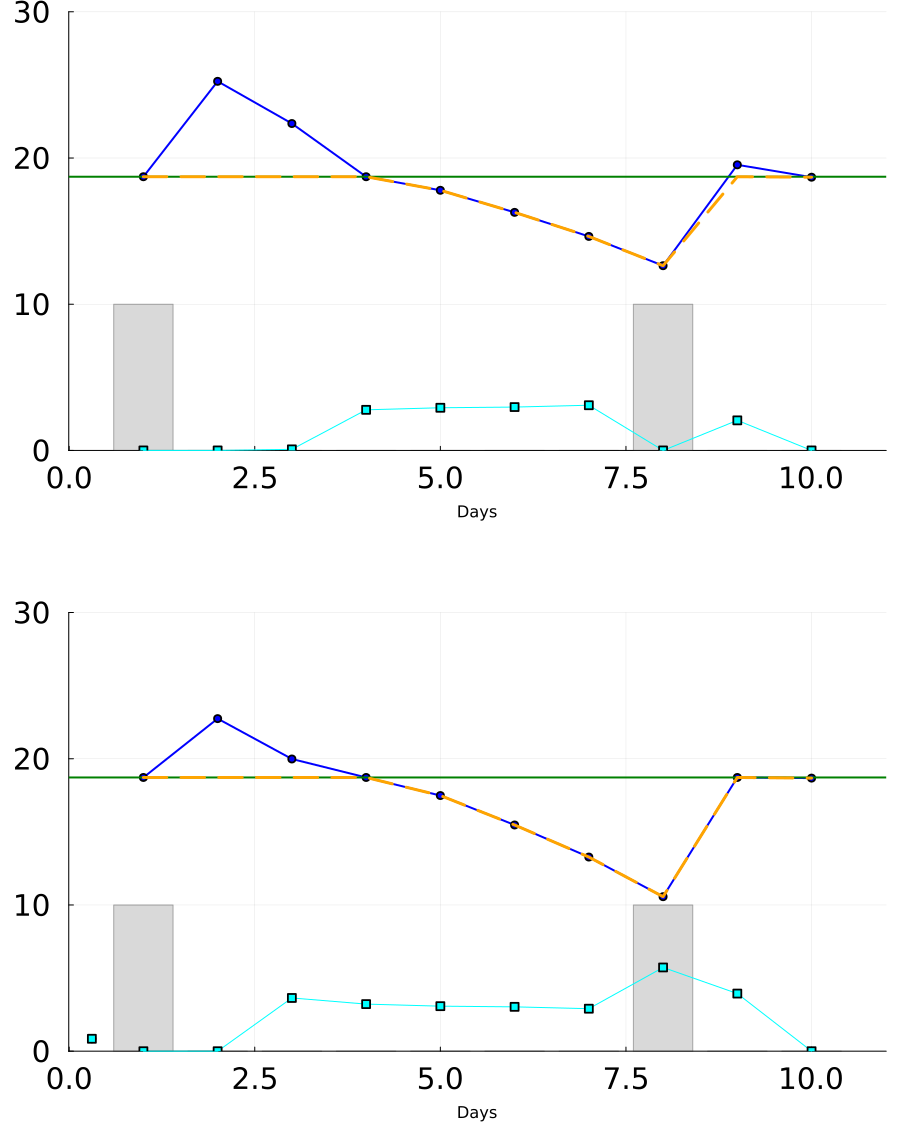}
    \caption{Irrigation profile for two fields with asymmetric slope topographies ($\alpha_1 = 0^{\circ}$, $\alpha_2 = 30^{\circ}$).}
    \label{fig:irrigation2}
\end{figure}

The second experiment introduces distinct environmental asymmetries by adjusting sector slopes ($\alpha_1 = 0^{\circ}$, $\alpha_2 = 30^{\circ}$). The resulting distinct execution dynamics are visualized in Fig.~\ref{fig:irrigation2}. The structural variations cause high runoff losses within the steep zone, prompting localized irrigation triggers on day 2, a full day ahead of the flat field sector. On day 8, irrigation occurs despite natural precipitation, indicating that the system accounts for cumulative water loss from slope runoff. The coordinator adjusts $X_{\min, j, i}$ dynamically to minimize aggregate deviations while satisfying the absolute global shared threshold constraint.

\begin{figure}[htbp]
    \centering
    \includegraphics[width=0.45\textwidth, keepaspectratio]{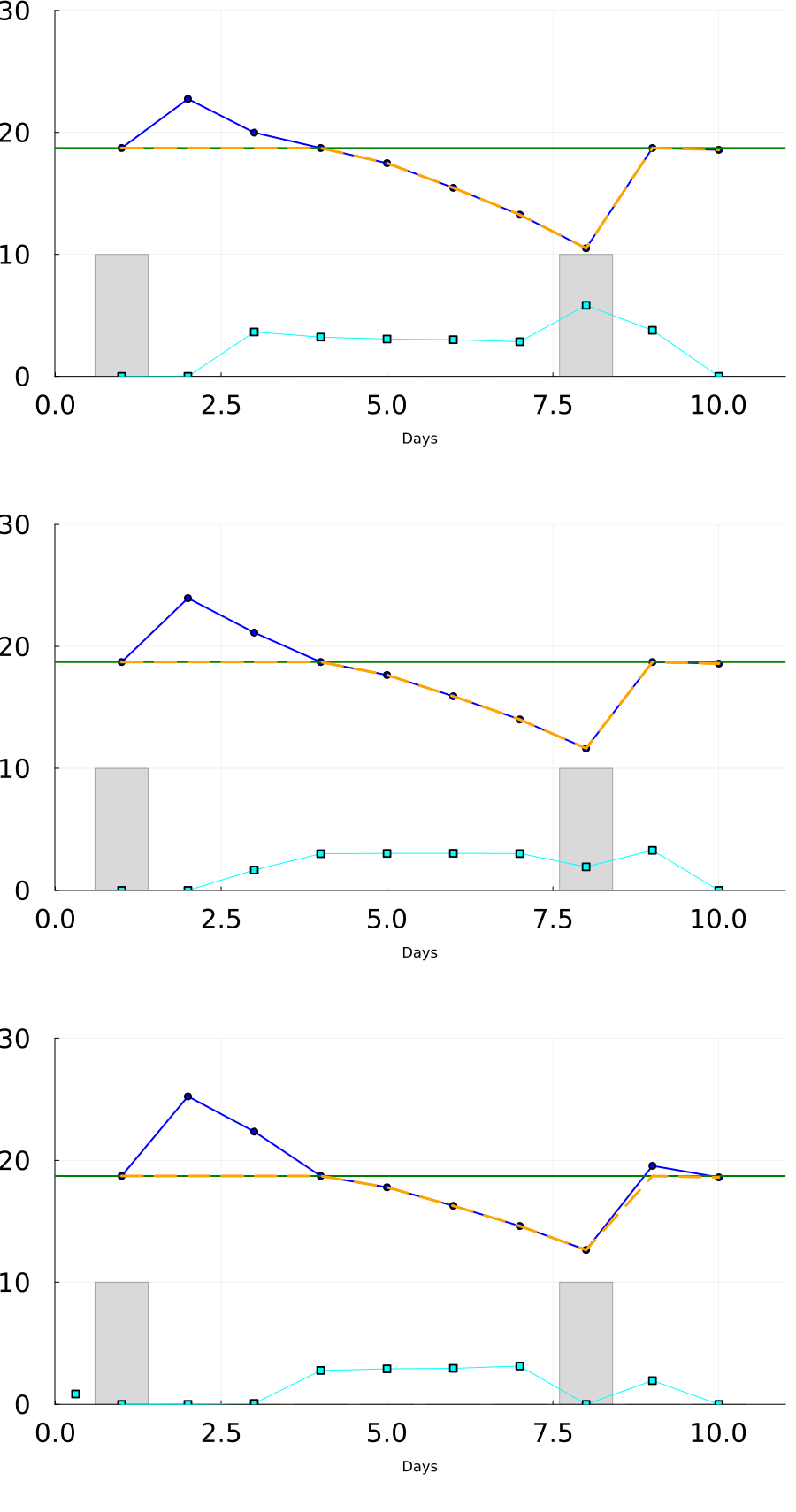}
    \caption{Scalability performance validation for three sectors with varying slope topographies ($\alpha = [30^{\circ}, 15^{\circ}, 0^{\circ}]$).}
    \label{fig:irrigation3}
\end{figure}

To evaluate scalability, a three-field scenario with slopes of $30^{\circ}$, $15^{\circ}$, and $0^{\circ}$ was evaluated under a scaled capacity limit $C_{\text{total}} = 9$~mm, as shown in Fig.~\ref{fig:irrigation3}.

Runoff losses map strictly to topographical variations: steeper plots show rapid depletion rates and initiate early water calls (day 2), whereas flatter sectors preserve moisture through day 3. This confirms that the bilevel structure balances competitive demands across an arbitrary number of points under constrained limits.


\section{Conclusion}
We presented a bilevel optimization framework that transforms smart irrigation from a rigid, point-based control problem into an adaptive, hierarchical resource management system. By promoting the soil moisture safety threshold to a strategic decision variable matrix $X_{\min, j, i}$, the model enables a  coordinator to dynamically adjust local moisture targets in response to global water scarcity and varying field topographies. Numerical results, obtained via the JuMP/IPOPT ecosystem on the Deucalion supercomputer, demonstrate that this approach ensures crop survival through controlled deficit irrigation while maximizing overall water-use efficiency. This framework provides a robust foundation for autonomous decision-support systems in large-scale, sustainable precision agriculture.

\section*{Acknowledgment}

The authors acknowledge financial support by national funds through the FCT/MCTES (PIDDAC), under the Associate Laboratory Advanced Production and Intelligent Systems -- ARISE LA/P/0112/2020 (DOI: 10.54499/LA/P/0112/2020) and the Base Funding (UIDB/00147/2020) and Programmatic Funding (UIDP/00147/2020) of the R\&D Unit Center for Systems and Technologies -- SYSTEC, as well as the projects NORTE2030-FEDER-02700700 and NORTE2030-FEDER-02719900.



\end{document}